\documentclass[12pt]{article}

\usepackage{fullpage}
\usepackage{amsmath}
\usepackage{amssymb}
\usepackage{amsthm}
\usepackage{graphicx}

\newtheorem{theorem}{Theorem}
\newtheorem{proposition}[theorem]{Proposition}

\newtheorem{definition}[theorem]{Definition}

\theoremstyle{remark}
\newtheorem{example}[theorem]{Example}

\numberwithin{theorem}{section}
\numberwithin{equation}{section}

\newcommand{\R}{{\mathbb R}}
\newcommand{\Rb}{\overline{\mathbb R}}
\newcommand{\E}{{\mathbf E}}
\renewcommand{\P}{{\mathbf P}}

\newcommand{\sA}{\mathcal{A}}
\newcommand{\sAR}{\mathcal{A}^{\mathrm{r}}}

\newcommand{\sP}{\mathcal{P}}
\renewcommand{\AA}{\mathbb{A}}
\newcommand{\AAR}{\mathbb{A}^{\mathrm{r}}}

\newcommand{\FF}{\mathbb{F}}
\newcommand{\sF}{\mathcal{F}}

\newcommand{\salg}{\mathfrak{F}}
\newcommand{\Lbase}{L}
\newcommand{\Ldi}[1][d]{\Lbase_#1^\infty}
\newcommand{\cone}{\mathbb{C}}
\newcommand{\risk}{\varrho}
\newcommand{\depth}[1][\alpha]{\mathrm{D}^#1}
\newcommand{\sriskbase}{s}

\newcommand{\driskbase}{\mathbb{D}}

\newcommand{\srisk}[1][\alpha]{\sriskbase_#1}
\newcommand{\sriskR}[1][\alpha]{\sriskbase_#1^{\mathrm{r}}}

\newcommand{\sriskn}[1][{1/n}]{\sriskbase_{#1}}
\newcommand{\drisk}[1][\alpha]{\driskbase^#1}
\newcommand{\arisk}[1][\alpha]{\sA_#1^\driskbase}
\renewcommand{\arisk}[1][\alpha]{\sA_#1}
\newcommand{\ariskR}[1][\alpha]{\sA_#1^{\mathrm{r}}}

\newcommand{\neutral}{\mathbf{e}}
\newcommand{\ftr}{\mathfrak{f}}
\newcommand{\Ind}{\mathbf{1}}

\newcommand{\VaR}[1][\alpha]{\mathrm{V@R}_{#1}}

\newcommand{\ES}[1][\alpha]{\mathrm{ES}_{#1}}
\newcommand{\EM}[1][{1/n}]{\mathrm{EM}_{#1}}
\newcommand{\HD}{\mathrm{HD}}
\newcommand{\ZD}{\mathrm{ZD}}
\newcommand{\CD}{\mathrm{CD}}
\newcommand{\essinf}{\mathrm{essinf}}
\newcommand{\Prob}[1]{\mathbf{P}\{#1\}}

\newlength{\querylen}
\usepackage{fancybox}

\newcommand{\pa}{\odot}
\newcommand{\ps}{\oplus}
\newcommand{\pmi}{\oplus}
\newcommand{\leqK}{\leq_{\scriptscriptstyle K}}
\newcommand{\coneK}{\cone_K}
\newcommand{\coneKR}{\cone_K^{\mathrm{r}}}
\newcommand{\KR}{K^{\mathrm{r}}}
\renewcommand{\top}{\mathbb{T}}

\begin{document}
\bibliographystyle{plain}

\title{Multivariate risks and depth-trimmed regions}

\author{Ignacio Cascos\footnote{Supported by the Spanish Ministry of
    Science and Technology under grant MTM2005-02254.}\\
  {\itshape\normalsize Department of Statistics, Universidad Carlos
    III de Madrid,}\\{\itshape\normalsize
    Av. Universidad 30, E-28911 Legan\'es (Madrid), Spain} \\[5mm]
  Ilya Molchanov\footnote{Supported by Swiss National
    Science Foundation.}\\ {\itshape\normalsize Department of
    Mathematical Statistics and Actuarial
    Science,}\\{\itshape\normalsize University of Berne, Sidlerstr. 5,
    CH-3012 Berne, Switzerland}\footnote{Corresponding author: Ilya
    Molchanov, Department of Mathematical Statistics and Actuarial
    Science, University of Berne, Sidlerstr. 5, CH-3012 Berne,
    Switzerland, e-mail: ilya.molchanov@stat.unibe.ch, Fax: 0041 31
    631 3870}}

\date{}
\maketitle


\begin{abstract}
  We describe a general framework for measuring risks, where the risk
  measure takes values in an abstract cone. It is shown that this
  approach naturally includes the classical risk measures and
  set-valued risk measures and yields a natural definition of
  vector-valued risk measures. Several main constructions of risk
  measures are described in this axiomatic framework.

  It is shown that the concept of depth-trimmed (or central) regions
  from the multivariate statistics is closely related to the
  definition of risk measures. In particular, the halfspace trimming
  corresponds to the Value-at-Risk, while the zonoid trimming yields
  the expected shortfall. In the abstract framework, it is shown how
  to establish a both-ways correspondence between risk measures and
  depth-trimmed regions. It is also demonstrated how the lattice
  structure of the space of risk values influences this relationship.
\end{abstract}



\noindent AMS Classification {91B30 \and 91B82 \and 60D05 \and 62H99}

\section{Introduction}
\label{sec:introduction}

Risk measures are widely used in financial engineering to assess the
risk of investments and to optimise the capital allocation. The modern
theory of coherent risk measures \cite{art:del:eber:99,delb02} aims to
derive properties of risk measures from several basic axioms:
translation-invariance, monotonicity, homogeneity and convexity. The
risk measures are mostly considered in the univariate case, i.e. it is
assumed that all assets have been transferred to their monetary
values. The quantile-based risk measures gain a particular importance
in the form of so-called spectral risk measures that are weighted
integrals of the quantile function, see \cite{acer02}.

When assessing risks of multivariate portfolios, the situation becomes
more complicated. The quantile function is not a numerical function
any more, and it is not possible to represent all portfolios as
functions of a uniform random variable.  The simplest approach to
assess the risk of a multivariate portfolio is to aggregate the
individual assets using their cash equivalents and then assess the
risk of the combined univariate portfolio. Then all portfolios with
identically distributed monetary equivalents would have identical
risks.

Several recent papers suggest various alternative ways of measuring
risks for multivariate portfolios without taking their monetary
equivalents.  The multivariate analogue of the Value-at-Risk discussed
in \cite{emb:puc06} is based on set-valued quantiles of the
multivariate cumulative distribution function.  A construction of
real-valued multivariate risk measures based on combining univariate
risks from transformed portfolios is described in \cite{bur:rues06}.
Multivariate coherent risk measures have been studied in
\cite{jouin:med:touz04} following the techniques from \cite{delb02}
based on the duality representations. The risk measures considered in
\cite{jouin:med:touz04} are actually set-valued and the preference
order corresponds to the ordering of sets by inclusion. It is
interesting to note that this order has the same meaning for risk, but
formally is the exactly opposite to the ordering of univariate risks
from \cite{art:del:eber:99}. Set-valued risk measures have been also
studied in \cite{ham06}. 

Because of this reason and in order to unify several existing
definitions we decided to consider risk measures as maps that have
values in a certain partially ordered cone, which may be, e.g., the
real line or the Euclidean space or the family of convex sets in the
Euclidean space. We single out the basic properties of so defined risk
measures and then describe the main technical constructions that make
it possible to produce new risk measures from the existing ones while
respecting their properties, e.g. the homogeneity or coherence. It is
not always assumed that the risk measures are coherent. Note that risk
measures with values in a partially ordered cone have been considered
in \cite{jas:kuec01}, where however it was assumed that this cone is
embeddable into a linear space. This is not the case for set-valued
risk measures which are also covered by the current work. These
set-valued measures can be used to produce vector-valued or
real-valued risk measures for multivariate portfolios.


In comparison with the studies of multivariate risk measures, the
multivariate statistical theory has an impressive toolbox suitable to
handle random vectors. We show that the multivariate setting for the
risk measures has a number of common features with the concept of
central (or depth-trimmed) regions well known in multivariate
statistics \cite{zuo:ser00a,zuo:ser00b}. They associate a random
vector with a set formed by the points in space located near to the
``central value'' of this random vector.  The risk measure is
generated by considering all translations of a random vector that
bring its central region to the positive (acceptable) part of the
space. In other words, the risks is determined by the relative
location of the central region comparing to the acceptable or
completely non-acceptable risks. Note that in the multivariate setting
the sets of acceptable and non-acceptable risk values are no longer
complementary, as they are in the univariate setting. Estimation
methods for depth trimmed regions then may be utilised to come up with
estimators for multivariate risk measures.  Despite the fact that the
definition of central regions (and indeed the name also) treats all
directions in the same way, it is possible to establish a two-way link
between depth-trimmed regions and risk measures.

\medskip

The paper is organised as follows.  Section~\ref{sec:risk} introduces
the main concept of a risk measure with values in an abstract cone. As
special cases one obtains the classical risk measures
\cite{art:del:eber:99}, set-valued risk measures of
\cite{jouin:med:touz04} and vector-valued risk measures.  A crucial
concept here is the function that assigns risks to deterministic
outcomes and controls changes of the risk if a deterministic amount is
being added to a portfolio. The partial order relation on the space of
risks makes it possible to consider it as a lattice.

The acceptance cone constitutes a subset of acceptable values for the
risk measure, while the acceptance set is the family of random vectors
whose risks belong to the acceptance cone.
Section~\ref{sec:acceptance-sets} discusses the main properties of the
acceptance set and the acceptance cone. We single out conditions that
make it possible to retrieve the risk measure from the acceptance set
it generates. This self-consistency condition can be traced to some
facts from the morphological theory of lattices \cite{hei94}.

Section~\ref{sec:constr-risk-meas} describes several ways to
construct new risk measures: minimisation, re-centring,
homogenisation, worst conditioning and transformations of risks.
In particular, the worst conditioning is a generic construction
that yields the expected shortfall if applied to the expectation.
It is shown that by transforming risks it is possible to produce
vector-valued risk measures from set-valued risk measures. This
construction can be applied, for instance, to the set-valued risk
measures from \cite{jouin:med:touz04}.

The definition of depth-trimmed regions and their essential properties
in view of relationships to risk measures are given in
Section~\ref{sec:depth}. In particular, the well-known halfspace
trimmed regions \cite{mas:theod94,rous:rut99} correspond to the
Value-at-Risk and the zonoid trimming \cite{mos02} produces the
expected shortfall. This analogy goes much further and leads to a
systematic construction of a risk measure from a family of
depth-trimmed regions in Section~\ref{sec:accept-sets-gere}. The main
idea here is to map the depth-trimmed region of a random vector into
the risk space using the function that assigns risks to deterministic
outcomes and then consider all translations of the image (of the
depth-trimmed region) that place it inside the acceptance cone.
Examples of basic risk measures obtained this way are described in
Section~\ref{sec:basic-risk-measures}. It is shown in
Section~\ref{sec:depth-trimm-regi} that the correspondence between
risk measures and depth-trimmed regions goes both ways, i.e.  it is
possible to construct a family of depth-trimmed regions from a risk
measure, so that, under some conditions, the initial risk measure is
recoverable from the obtained family of depth-trimmed regions.

Finally, Section~\ref{sec:duality-results} deals with dual
representation of coherent risk measures and depth-trimmed regions
using families of measures, in a way similar to the well-known
approach \cite{delb02} for real-valued coherent risk measures. In
particular we show that all coherent vector-valued risk measures are
marginalised, i.e. can be represented as the vector of risk measures
for the marginals. This fact confirms the idea that set-valued risk
measures are a natural tool for multivariate portfolios if one is
interested in non-trivial coherent risk measures.

\section{Risk measures in abstract cones}
\label{sec:risk}

A risky portfolio is modelled as an essentially bounded random vector
$X$ that represents a financial gain.  Let $\Ldi$ denote the set of
all essentially bounded $d$-dimensional \emph{random vectors} on the
probability space $(\Omega,\salg,\P)$. In order to combine several
definitions of risk measures, it is sensible to regard them as
functionals on $\Ldi$ with values in a partially ordered \emph{convex
  cone} $\cone$.

\begin{definition}[Semigroup and convex cone]
  \label{def:cc}
  An \emph{abelian topological semigroup} is a topological space
  $\cone$ equipped with a commutative and associative continuous
  binary operation $\ps$. It is assumed that $\cone$ possesses the
  \emph{neutral element} $\neutral$ satisfying $x\ps \neutral=x$ for
  all $x\in\cone$.  The semigroup $\cone$ is a \emph{convex cone} if
  it is also equipped with a continuous operation $(x,t)\mapsto t\pa
  x$ of multiplication by positive scalars $t>0$ for $x\in\cone$ so
  that $1\pa x=x$ for all $x\in\cone$, $t\pa \neutral=\neutral$ for
  all $t>0$, and the following conditions are satisfied
  \begin{align*}
    t\pa (x\ps y)&=t\pa x\ps t\pa y, \quad t>0,\; x,y\in\cone,\\
    t\pa(s\pa x)&=(ts)\pa x,\quad t,s>0,\; x\in\cone.
  \end{align*}
\end{definition}

Assume throughout that $\cone$ is endowed with a \emph{partial order}
$\preceq$ that is compatible with the (commutative) addition operation
and multiplication by scalars, i.e.  $x\preceq y$ implies that $x\ps
z\preceq y\ps z$ for all $z$ and $t\pa x\preceq t\pa y$ for all $t>0$.
Furthermore, assume that $\cone$ with the order $\preceq$ is a
\emph{complete lattice}, i.e. every set has supremum and infimum,
which are denoted by $\vee$ and $\wedge$ respectively.  Since this
partial order may differ from the conventional order for real numbers,
we retain the notation supremum and infimum (also min and max) for the
conventional order on the real line, while $\vee$ and $\wedge$ denote
the supremum and infimum in $\cone$.  The \emph{top element} of
$\cone$ is denoted by $\top$.  It is assumed that the top element is
absorbing, i.e. $\top\ps a=\top$ for all $a\in\cone$.

Note that the cone $\cone$ is not necessarily embeddable in a linear
space, since the addition operation does not necessarily obey the
cancellation law and the second distributivity law $t\pa x\ps s\pa
x=(t+s)\pa x$ is not imposed, see \cite{dav:mol:zuy05} for a
discussion of algebraic properties of convex cones. Accordingly, it is
not possible to view $\cone$ as a partially ordered linear space. This
situation is typical if $\cone$ is the family of convex sets in the
Euclidean space $\R^d$ and the additive operation is the \emph{closed
  Minkowski addition}, i.e. the sum of $A\pmi B$ of two sets is the
topological closure of $\{x+y:\; x\in A,\; y\in B\}$.  Note that the
Minkowski sum of two non-compact closed sets is not necessarily
closed.  The multiplication by positive numbers is given by
$tA=\{tx:\: x\in A\}$, i.e. the usual dilation of $A$ by $t>0$ and we
simply write $x+A$ instead of $\{x\}\pmi A$.

We retain the usual $+$ and multiplication signs for operations with
real numbers and vectors in $\R^d$. For convenience, letters $x,y,z$
with or without subscripts stand for points in $\R^d$, letters $t,s$
represent real numbers, letters $a,b$ denote elements of $\cone$,
letters $X,Y$ are used for random variables or random vectors, and
$A,B,F,K$ are subsets of $\R^d$.

A \emph{proper Euclidean convex cone} $K$ is a strict subset of $\R^d$
such that $\R^d_+\subseteq K$, $K$ does not contain any line, and
$x+y\in K$, $tx\in K$ for all $x,y\in K$ and $t>0$. In the univariate
case ($d=1$) the only possibility is $K=[0,\infty)$.

\begin{definition}[Order in $\R^d$]
  \label{def:rd-order}
  Let $K$ be a proper Euclidean convex cone. For $x,y\in\R^d$, we
  denote $x\leqK y$ if and only if $y-x\in K$.
\end{definition}

From the economical viewpoint, this ordering would correspond, e.g. to
exchanges of various currencies, cf. \cite{jouin:med:touz04,kab99}.

\medskip

A risk measure is a functional on $\Ldi$ with values in $\cone$.  As
the first step of its proper definition, one should specify how this
functional acts on degenerate random variables, i.e. on the space
$\R^d$, which is naturally embedded in $\Ldi$. This action is defined
by a function $f:\R^d\mapsto\cone$, which is interpreted as the risk
associated with the degenerate random variable $X=x$ a.s. Assume that
$f(0)=\neutral$, $f$ is linear, i.e.
\begin{equation}
  \label{eq:f}
  f(x)\ps f(y)=f(x+y)
\end{equation}
for all $x,y\in\R^d$, and non-decreasing, i.e. $f(y)\preceq f(x)$ if
$y\leqK x$. The mapping $f$ is a linear positive map between partially
ordered linear spaces: $\R^d$ with the $\leqK$ order and the space
$\FF=\{f(y):\; y\in\R^d\}$ with the order inherited from $\cone$.
Condition (\ref{eq:f}) implies that $f(x)\neq\top$ for all $x$.
Indeed, if $f(x)=\top$, then $f(x+y)=\top$ for all $y$, so that $f$
identically equals $\top$ contrary to the fact that $f(0)=\neutral$.

The following definition specifies the desirable properties of
risk measures.

\begin{definition}[Risk measure]
  \label{def:1}
  A functional $\risk:\Ldi\mapsto\cone$ is called a \emph{risk
    measure} associated with $f$ if $\risk(X)=f(x)$ in case $X=x$ a.s.
  and the following conditions hold
  \begin{description}
  \item[\bf R1] $f(y)\ps\risk(X)=\risk(X+y)$ for all $y\in\R^d$;
  \item[\bf R2] $\risk(Y)\preceq \risk(X)$ whenever $Y\leqK X$ a.s.;
  \end{description}
  is called a \emph{homogeneous} risk measure if also
  \begin{description}
  \item[\bf R3] $\risk(tX)=t\pa \risk(X)$ for all $t>0$ and $X\in\Ldi$;
  \end{description}
  and a \emph{coherent risk measure} if additionally
  \begin{description}
  \item[\bf R4] $\risk(X)\ps \risk(Y)\preceq\risk(X+Y)$ for all $X,Y\in\Ldi$.
  \end{description}
\end{definition}


Condition \textbf{R2} means that $\risk$ is a lattice morphism between
$\Ldi$ with the partial order generated by $\leqK$ and $\cone$.  It is
also possible to consider not necessarily homogeneous risk measures
that satisfy the assumption
\begin{equation}
  \label{eq:crm}
  t\pa\risk(X)\ps(1-t)\pa\risk(Y)\preceq\risk(tX+(1-t)Y)
\end{equation}
for all $t\in[0,1]$, which are traditionally called \emph{convex}
\cite{foel:sch02c} (despite the fact that the inequality in our
setting actually means that $\risk$ is concave).

Note that the multiplication by numbers in $\cone$ is not needed if
\textbf{R3} is not considered. In this case one can only require that
$\cone$ is a partially ordered abelian semigroup.  Furthermore,
Definition~\ref{def:1} can be formulated for any partially ordered
cone $\cone$ (not necessarily a complete lattice) and any partial
ordering on $\R^d$.

Since $\neutral\neq\top$, the condition $f(0)=\neutral$ together with
\textbf{R2} imply that $\risk(X)$ never takes the value $\top$. This
corresponds to the requirement that conventional risk measures do not
take the value $-\infty$, see \cite{delb02}. Indeed, if
$\risk(X)=\top$, then $f(a)=\risk(a)=\top$ for $a$ being an upper
bound for $X$.

The use of function $f$ in Definition~\ref{def:1} is twofold.  It
determines risks of deterministic portfolios and also controls how the
risk of $X$ changes if a deterministic quantity is added to the
portfolio $X$.  The second task can be also delegated to another
function $g:\R^d\mapsto\cone$, so that \textbf{R1} becomes $g(y)\ps
\risk(X)=\risk(X+y)$ and $g(y)\ps g(-y)=\neutral$ for all $y\in\R^d$.
It is easy to show that $f$ and $g$ coincide if and only if
$f(0)=\neutral$.


\begin{example}[Set-valued risk measures]
  \label{ex:jouini}
  Consider the family of closed convex sets in $\R^d$ partially
  ordered by inclusion with the addition defined as the closed
  Minkowski sum and the conventional dilation by positive numbers.
  Define $f(x)=\{y\in\R^d:\; -x\leqK y\}=-x+K$, where $K$ is a proper
  Euclidean cone from Definition~\ref{def:rd-order}.  In particular,
  the fact that $\risk(X)\supset K$ means that $X$ has a negative
  risk.  In this case Definition~\ref{def:1} turns into
  \cite[Def.~2.1]{jouin:med:touz04}. Since $f(0)=K$ has to be the
  neutral element, the relevant cone $\cone$ should consist of all
  closed convex sets $F\subseteq\R^d$ such that the closed Minkowski
  sum $F\pmi K$ coincides with $F$. This important family of sets will
  be denoted by $\coneK$.

  Let us show that $\coneK$ is a complete lattice. Consider any family
  of sets $\{A_i:\;i\in I\}\subseteq\coneK$. Then $F=\bigvee_{i\in I}
  A_i$ is the smallest convex set that contains all the $A_i$s, i.e.
  $F$ is the closure of the convex hull of the union of these sets.
  Since $F$ is closed convex and $F=\bigvee (A_i\pmi K)=K\pmi F$, we
  have $F\in\coneK$. Furthermore, $M=\bigwedge_{i\in I} A_i$ is given
  by $M=\bigcap_{i\in I} A_i$. The set $M$ is closed convex and also
  belongs to $\coneK$, since
  \begin{displaymath}
    M=\bigcap_{i\in I} A_i=\bigcap_{i\in I} (A_i\pmi K)
    \supseteq K\pmi \bigcap_{i\in I} A_i
    \supseteq M\,,
  \end{displaymath}
  because $K$ contains the origin.
\end{example}

\begin{example}[Univariate risk measures]
  \label{ex:artzner}
  The classical definition of real-valued coherent risk measures from
  Artzner \emph{et al.}~\cite{art:del:eber:99} can be recovered from
  the setting of Example~\ref{ex:jouini} for $d=1$ and
  $\risk(X)=[\rho(X),\infty)$, where $\rho(X)$ is the risk measure of
  $X$ as in \cite{art:del:eber:99}. An alternative approach is to let
  $\cone$ be the extended real line $\Rb=[-\infty,\infty]$ with the
  \emph{reversed} order and conventional addition and multiplication
  operations. In this case $f(x)=-x$.  We will briefly recall three
  univariate risk measures: the value at risk, which is the most
  widely used risk measure, and two coherent risk measures, the
  expected shortfall and the expected minimum.

  The \emph{value at risk} is defined as the amount of extra capital
  that a firm needs in order to reduce the probability of going
  bankrupt to a fixed threshold $\alpha$. It is the opposite of the
  $\alpha$-quantile of a random variable $X$, i.e.
  \begin{displaymath}
    \VaR(X)=-\inf\{x:\; \Prob{X\leq x}>\alpha\}=-F^{-1}_X(\alpha)\,,
  \end{displaymath}
  where $F_X$ is the cumulative distribution function of $X$.
  It can be shown that the value at risk is a homogeneous risk
  measure, but not a coherent one. It satisfies properties
  \textbf{R1}, \textbf{R2} and \textbf{R3}, but not necessarily
  \textbf{R4}.

  The \emph{expected shortfall} is a coherent risk measure defined as
  \begin{displaymath}
    \ES(X)= -\,\frac{1}{\alpha}\int_0^\alpha F_X^{-1}(t){\rm d}t\,,
  \end{displaymath}
  where $\alpha\in(0,1]$.

  The \emph{expected minimum} is another coherent risk measure defined
  as
  \begin{displaymath}
    \EM(X)= - \E\min\{X_1,X_2,\dots,X_n\},
  \end{displaymath}
  where $X_1,X_2,\dots,X_n$ are independent copies of $X$. The
  expected minimum belongs to the family of weighted V@Rs and is
  called Alpha-V@R in \cite{cher:mad06c}.
\end{example}


In the following we often consider the Euclidean space $\R^d$ extended
by adding to it the top and bottom elements at the infinity, so that
the space then becomes a complete lattice. In order to simplify
notation we retain notation $\R$ and $\R^d$ for such extended
spaces.

\begin{example}[Marginalised multivariate vector-valued risk measures]
  \label{ex:vector}
  Let $\cone$ be $\R^d$ with the usual addition, multiplication by
  positive numbers and the reversed coordinatewise order, i.e.
  $a\preceq b$ if $b\leqK a$ with $K=\R_+^d$. Given a $d$-dimensional
  random vector $X=(X_1,\dots,X_d)$, any of the aforementioned
  univariate risk measures $\rho$ yields a risk measure
  $\risk(X)=(\rho(X_1),\dots,\rho(X_d))$ with values in $\R^d$. In
  this case $f(x)=-x$.
\end{example}

\section{Acceptance cones and acceptance sets}
\label{sec:acceptance-sets}

The concept of an acceptance set is the dual one to the risk measure,
see \cite{art:del:eber:99,foel:sch02c,jouin:med:touz04}.  The main
idea is that a portfolio $X$ is acceptable if $\risk(X)$ belongs to a
certain subcone $\AA\subset\cone$ called the \emph{acceptance cone}.
The classical setting (see Example~\ref{ex:artzner}) corresponds to
$\cone=\R$ with the reversed order and $\AA=(-\infty,0]$.  Every
acceptance cone $\AA$ is upper with respect to $\preceq$, i.e.  if
$a\preceq b$ and $a\in\AA$, then $b\in\AA$. We also assume that
\begin{equation}
  \label{f-acceptance}
  \{a\in\cone:\; \neutral\preceq a\}=\AA\,,
\end{equation}
i.e. a deterministic portfolio $x$ is acceptable if and only if $0\leqK
x$.


Given the risk measure $\risk$, the set $\sA\subset\Ldi$ of acceptable
portfolios (called the \emph{acceptance set}) is given by
\begin{displaymath}
  \sA=\{X\in\Ldi:\; \risk(X)\in\AA\}
  =\{X\in\Ldi:\; \neutral\preceq \risk(X)\}\,.
\end{displaymath}
If $\risk$ is coherent, then $\sA$ is a cone in $\Ldi$.  It follows
from \textbf{R1} that
\begin{displaymath}
  \{y:\; X-y\in\sA\}=\{y:\; \risk(X-y)\in\AA\}
  =\{y:\; \risk(X)\ps f(-y)\in\AA\}\,.
\end{displaymath}
The $f$-image of the set in the right-hand side is
\begin{align*}
  \risk_\sA(X)&=\{f(y):\; y\in\R^d,\; \risk(X)\ps f(-y)\in\AA\}\\
  &=\{a\in\FF:\; \risk(X)\in \AA\ps a\}\\
  &=\{a\in\FF:\; a\preceq \risk(X)\}\,.
\end{align*}
Indeed, since the family $\FF$ of values of $f$ is a linear space,
$\AA\ps a=\{b\ps a:\; b\in\AA\}$ coincides with the set
$\{b\in\cone:\; a\preceq b\}$ for any $a\in\FF$.

Note that $\risk_\sA(X)$ is not necessarily an element of $\cone$,
since it may consist of several elements of $\cone$.  For instance, in
Example~\ref{ex:artzner} (with $\cone=\R$), $\risk_\sA(X)$ is the set
$[\rho(X),\infty)$, while the risk of $X$ is a real number.  In this
case, one can retrieve the risk of $X$ by taking the infimum of all
members of $\risk_\sA(X)$.  This minimum corresponds to the
$\vee$-operation in $\R$ with the reversed order.  The following easy
observation generalises the well-known relationship between risk
measures and acceptance sets \cite{art:del:eber:99,delb02}.

\begin{proposition}
  \label{prop:r-a}
  If $\FF$ is sup-generating~(see \cite[p.~28]{hei94}), i.e.
  \begin{equation}
    \label{eq:sg}
    b=\bigvee\{a\in\FF:\; a\preceq b\}\quad\textrm{for all}\;\;
    b\in\cone\,,
  \end{equation}
  then
  \begin{displaymath}
    \risk(X)=\bigvee \risk_\sA(X)\,.
  \end{displaymath}
\end{proposition}

In the multivariate case one often needs the concept of the rejection
cone $\AAR=\{a\in\cone:\; a\preceq\neutral\}$ and the rejection set
\begin{displaymath}
  \sAR= \{X\in\Ldi:\; \risk(X)\preceq \neutral\}\,.
\end{displaymath}
While $\AAR$ is a subcone of $\cone$, the set $\sAR$ is not
necessarily convex even if $\risk$ is coherent. Indeed, if
$X,Y\in\sAR$, then $\risk(X)+\risk(Y)\preceq\neutral$, while \textbf{R4}
no longer suffices to deduce that $\risk(X+Y)\preceq\neutral$.

\begin{example}[Set-valued risk measures]
  \label{ex:mvar-cones}
  Let $\coneK$ be the cone of convex closed sets described in
  Example~\ref{ex:jouini} and $f(x)=-x+K$, so that $\FF=\{y+K:\;
  y\in\R^d\}$. If $\AA=\{A\in\coneK:\; K\subseteq A\}$, then $\FF$ is
  sup-generating, since for any $F\in\coneK$ we have
  \begin{displaymath}
    \bigvee\{a\in\FF:\; a\preceq b\}
    =\bigcup\{y+K:\; y\in\R^d,\; (y+K)\subseteq F\}
    =F\,.
  \end{displaymath}
  As in \cite[Sec.~2.5]{jouin:med:touz04}, it is possible to choose
  another acceptance cone $\AA'$ which is richer than the cone $\AA$
  defined above. Furthermore, the sup-generating property
  (\ref{eq:sg}) corresponds to the self-consistency property from
  \cite[Property~3.4]{jouin:med:touz04}.
\end{example}

\begin{example}[Alternative construction of set-valued risk measures]
  \label{ex:ac-sval}
  There is also an alternative way to introduce set-valued risk
  measures. Let $\coneKR$ be the family of complements to the
  interiors of sets from $\coneK$, with the addition operation induced
  by one from $\coneK$, i.e. $F_1\ps F_2$ is the complement to the
  Minkowski sum of the complements to $F_1$ and $F_2$. The neutral
  element $\neutral=\KR$ is then the complement to the interior of
  $K$. If $\coneKR$ is equipped with the inclusion order, then the
  same arguments as in Example~\ref{ex:jouini} confirm that $\coneKR$
  is a complete lattice.

  If $f(x)=x+\KR$, $x\in\R^d$, then the corresponding family $\FF$ is
  inf-generating (see  \cite[p.~28]{hei94}), i.e.
  \begin{displaymath}
    b=\bigwedge\{a\in\FF:\; b\preceq a\}\quad\textrm{for all}\;\;
    b\in\cone\,.
  \end{displaymath}
  In this case
  \begin{equation}
  \label{eq:r-ar}
    \risk(X)=\bigwedge \risk_{\sAR}(X)\,,
  \end{equation}
  where $\risk_{\sAR}(X)$ is the $f$-image of all $y\in\R^d$ such that
  $X-y\in\sAR$. 
\end{example}

\begin{example}[Vector-valued risk measures from scalar portfolios]
  Consider a risk measure $\risk$ defined on $\Ldi[1]$ with values in
  $\cone=\R^2$ with the usual summation and multiplication by scalars
  and the reversed coordinatewise ordering, i.e. the reversed ordering
  generated by $K=\R_+^2$. Such a risk measure may be defined as a
  vector composed of several univariate risk measures from
  Example~\ref{ex:artzner}. In this case $f(x)=(-x,-x)$, so that $\FF$
  is the diagonal in $\R^2$, which is clearly not sup-generating.

  This example explains, by the way, why in the framework of
  \cite{jouin:med:touz04} only risk measures that do not increase the
  dimension of the portfolios have been studied.
\end{example}


\section{Constructions of risk measures}
\label{sec:constr-risk-meas}


\subsection{Minimisation}
\label{sec:minimisation}

Consider a family $\risk_i$, $i\in I$, of risk measures on the same
cone $\cone$, all associated with the same function $f$.
Then $\risk=\bigwedge_{i\in I}\risk_i$ is also a risk measure
associated with $f$. If all $\risk_i$ are coherent (resp.  homogeneous
or convex) the resulting risk measure is coherent (resp.  homogeneous
or convex).  The acceptance set associated with $\risk$ is the
intersection of the acceptance sets of the risk measures $\risk_i$,
$i\in I$.


\begin{example}[Minimisation of univariate risk measures]
  While it is not interesting to take minimum of, say, the expected
  shortfalls at different levels, it is possible to combine members
  from different families of univariate risk measures. For instance,
  if $n\geq 1$ and $\alpha\in(0,1]$, then $\max\{\EM(X),\ES(X)\}$ is a
  coherent risk measure associated with $f(x)=-x$. Note that the
  maximum of two risk measures correspond to the minimum in $\cone=\R$
  with the reversed order.
\end{example}

\subsection{Re-centring}
\label{sec:recentring}

All random vectors from $\Ldi$ can be naturally centred by subtracting
their expected values. This makes it possible to define a risk measure
on centred random vectors and then use \textbf{R1} to extend it onto
the whole $\Ldi$.  If $\risk$ is defined on the family of essentially
bounded random vectors with mean zero, then the re-centred risk
measure is given by
\begin{displaymath}
  \risk_{\mathrm{o}}(X)=\risk(X-\E X)\ps f(\E X)\,, \quad X\in \Ldi\,.
\end{displaymath}
If $\cone$ is $\R^d$ or a family of subsets of $\R^d$, we rely on the
canonical choice of the translation by setting
$\risk_{\mathrm{o}}(X)=\risk(X-\E X)-\E X$.


It should be noted that \textbf{R2} does not hold automatically for
re-centred risk measures and has to be checked every time the
re-centring is applied.

\subsection{Homogenisation}
\label{sec:homogenisation}

If $\risk$ satisfies \textbf{R1} and \textbf{R2}, it is possible to
construct a homogeneous risk measure from it by setting
\begin{equation}
  \label{eq:hc}
  \risk_{\mathrm{h}}(X)=\bigwedge_{t>0} \frac{1}{t}\pa\risk(tX)\,.
\end{equation}
Note that the infimum operation $\bigwedge$ in $\cone$
makes sense, since $\cone$ is a complete lattice. It is
easy to see that $\risk_{\mathrm{h}}$ satisfies \textbf{R3}.
Furthermore, it satisfies \textbf{R2} and \textbf{R1} if $f$ is
homogeneous. The latter is clearly the case if $f(x)=-x+K$,
$x\in\R^d$, for a proper cone $K$, see Example~\ref{ex:jouini}.

A similar construction produces a translation-invariant risk measure
from a general one by
\begin{equation}
  \label{eq:htc}
  \risk_{\mathrm{t}}(X)=\bigwedge_{z\in \R^d}
  \Big(\risk(X+z)\ps f(-z)\Big)\,.
\end{equation}
Both (\ref{eq:hc}) and (\ref{eq:htc}) applied together to a function
$\risk$ that satisfies \textbf{R2} and \textbf{R4} yield a coherent
risk measure.

\begin{example}
  \label{ex:homog}
  If $\cone$ is the real line with the reversed order and
  (\ref{eq:hc}) results in a non-trivial function, then $\risk(tX)\to
  0$ as $t\to 0$.  Similarly, a non-trivial result of (\ref{eq:htc})
  yields that $\risk(X+z)\to -\infty$ as $z\to\infty$.  For instance,
  these constructions produce trivial results if applied to the risk
  measure $\E (k-X)_+$ studied in \cite{jar02}.
\end{example}

\subsection{Worst conditioning}
\label{sec:worst-conditioning}

A single risk measure $\risk$ can be used to produce a family of risk
measures by taking the infimum of the risks associated to the random
vectors obtained after certain rearrangements of the underlying
probability measure. For each $\alpha\in(0,1]$ define
\begin{displaymath}
  \risk_\alpha(X)=\bigwedge_{\phi\in\Phi_\alpha} \risk(X_\phi)\,,
\end{displaymath}
where $X_\phi=X\circ\phi$ and $\Phi_\alpha$ is the family of
measurable mappings $\phi:\Omega\mapsto\Omega$ such that
$\P(\phi^{-1}(A))\leq\alpha^{-1}\P(A)$ for all $A\in\salg$. If
$X\in\Ldi$, then $X_\phi\in\Ldi$ for any $\alpha\in(0,1]$ and
$\phi\in\Phi_\alpha$.  It is possible to define the worst conditioning
alternatively as
\begin{displaymath}
  \risk_\alpha(X)=\bigwedge_{Y\in \sP_\alpha(X)} \risk(Y)\,,
\end{displaymath}
where $\sP_\alpha(X)$ is the family of all random vectors $Y$ such
that $\Prob{Y\in B}\leq \alpha^{-1}\Prob{X\in B}$ for all Borel
$B\subset\R^d$.

It is easy to show that $\risk_\alpha$ preserves any property that
$\risk$ satisfies from \textbf{R1}--\textbf{R4}. For instance, if
$Y\leqK X$ a.s., then $Y_\phi\leqK X_\phi$ a.s. for any
$\phi\in\Phi_\alpha$, so that
\begin{displaymath}
  \risk_\alpha(Y)=\bigwedge_{\phi\in\Phi_\alpha}
  \risk(Y_\phi)\preceq\bigwedge_{\phi\in\Phi_\alpha}
  \risk(X_\phi)=\risk_\alpha(X)
\end{displaymath}
whenever $\risk$ satisfies \textbf{R2}. If $X,Y\in\Ldi$ and $\risk$
satisfies \textbf{R4}, then
\begin{align*}
  \risk_\alpha(X+Y)&=\bigwedge_{\phi\in\Phi_\alpha}
  \risk\big((X+Y)_\phi\big)=\bigwedge_{\phi\in\Phi_\alpha}
  \risk(X_\phi+Y_\phi)\\
  &\succeq\bigwedge_{\phi\in\Phi_\alpha}
  \big(\risk(X_\phi)\ps \risk(Y_\phi)\big)
  \succeq\bigwedge_{\phi\in\Phi_\alpha}
  \risk(X_\phi)\ps \bigwedge_{\phi\in\Phi_\alpha}
  \risk(Y_\phi)=\risk_\alpha(X)\ps \risk_\alpha(Y)\,.
\end{align*}

Consider now the setting of univariate risk measures from
Example~\ref{ex:artzner}, i.e. $X$ is a random variable from $\Ldi[1]$
and $\cone$ is the real line with the reversed order.  The simplest
coherent risk measure is the opposite of the expectation of a random
variable. In fact, this risk measure appears from the expected
shortfall when $\alpha=1$, i.e.  $\ES[1](X)=-\E X$.  The worst
conditioning applied to the opposite of the expectation yields
\begin{align*}
  (-\E)_\alpha(X)&=\sup_{\phi\in\Phi_\alpha}\{-\E(X_\phi)\}
  =-\inf_{\phi\in\Phi_\alpha}\E(X_\phi)
  =-\inf_{\phi\in\Phi_\alpha}\int X(\phi(\omega))\P({\rm d}\omega)\\
  &=-\inf_{\phi\in\Phi_\alpha}\int X(\omega)\P\phi^{-1}({\rm
    d}\omega)=-\inf_{\phi\in\Phi_\alpha}\E_{\P\phi^{-1}}X\,,
\end{align*}
where $\E_{\P\phi^{-1}}$ denotes the expectation with respect to the
probability measure $\P\phi^{-1}$.



In general, $-\inf_{\phi\in\Phi_\alpha}\E_{\P\phi^{-1}}X\leq\ES(X)$
with the equality if $(\Omega,\salg,\P)$ is non-atomic.  Without loss
of generality assume that $\Omega=[0,1]$, $\P$ is the Lebesgue measure
restricted to $[0,1]$ and $X$ is increasing mapping from $[0,1]$ into
$\R$, which implies that $X(\omega)=F^{-1}_X(\omega)$ for all
$\omega\in[0,1]$, where $F_X$ is the cumulative distribution function
of $X$.  The infimum of $\E_{\P\phi^{-1}}X$ over all
$\phi\in\Phi_\alpha$ is achieved when $X\circ\phi$ takes the smaller
possible values with the highest possible probabilities, and thus it
is attained at $\phi'(\omega)=\alpha\omega$. We conclude
\begin{displaymath}
  (-\E)_\alpha(X)=-\int X(\alpha\omega){\rm d}\omega
  =-\frac{1}{\alpha}\int_0^\alpha F_X^{-1}(t){\rm d}t=\ES(X)\,,
\end{displaymath}
i.e. the expected shortfall appears by applying the worst conditioning
construction to the opposite of the expectation.

\begin{example}[Worst conditioning of the expected shortfall]
  Let us now apply the worst conditioning to the expected shortfall at
  level $\beta$,
  \begin{align*}
    \big(\ES[\beta]\big)_\alpha(X)
    &=\sup_{\phi_1\in\Phi_\alpha} \ES[\beta](X_{\phi_1})
    =\sup_{\phi_1\in\Phi_\alpha}\Big(-\inf_{\phi_2\in\Phi_\beta}
    \E_{\P\phi_1^{-1}}X_{\phi_2}\Big)\\
    &=-\inf_{\phi_1\in\Phi_\alpha,\;\phi_2\in\Phi_\beta}
    \E_{\P\phi_1^{-1}\phi_2^{-1}}X\,.
  \end{align*}
  Clearly $\phi_2\circ\phi_1\in\Phi_{\alpha\beta}$ and thus
  $\big(\ES[\beta]\big)_\alpha(X)\leq\ES[\alpha\beta](X)$.  If the
  probability space is non-atomic, all mappings from
  $\Phi_{\alpha\beta}$ can be written as the composition of a mapping
  from $\Phi_\alpha$ and a mapping from $\Phi_\beta$, so that
  $\big(\ES[\beta]\big)_\alpha(X)=\ES[\alpha\beta](X)$. One can say
  that the expected shortfall is stable under the worst conditioning.
\end{example}

\begin{example}[Worst conditioned $\VaR$]
  Let us finally apply the worst conditioning construction to the
  value at risk at level $\beta$ considered on a non-atomic
  probability space $\Omega=[0,1]$ with $\P$ being the Lebesgue
  measure. Without loss of generality assume that $X$ is increasing,
  so that $X(\omega)=F^{-1}_X(\omega)$.  The infimum below is attained
  at $\phi'(\omega)=\alpha\omega$ and since $X_{\phi'}$ is also
  increasing, we have $X_{\phi'}(\omega)=F^{-1}_{X_{\phi'}}(\omega)$.
  Thus
  \begin{displaymath}
    \big(\VaR[\beta]\big)_\alpha(X)
    =-\inf_{\phi\in\Phi_\alpha}F^{-1}_{X_\phi}(\beta)
    =-X_{\phi'}(\beta)=-X(\alpha\beta)=-F^{-1}_X(\alpha\beta)
    =\VaR[\alpha\beta](X)\,.
  \end{displaymath}
\end{example}

\subsection{Transformations of risks}
\label{sec:transf-risks-1}

Risk measures with values in a cone $\cone$ may be further transformed
by mapping $\cone$ into another cone $\cone'$ using a map $h$. The aim
may be to change the dimensionality (cf. \cite{jouin:med:touz04}) or
produce a vector-valued risk measure from a set-valued one.  The map
$h:\cone\mapsto\cone'$ that transforms any $\cone$-valued risk measure
$\risk$, into the $\cone'$-valued risk measure $h(\risk(\cdot))$, will
be called a \emph{risk transformation}. If $h$ respects the coherence
property of risk measures, it will be called a \emph{coherent map}.

Let us denote by $\preccurlyeq$ the partial order in $\cone'$ which we
assume to be compatible with the (commutative) addition operation and
multiplication by scalars. The additive operation on $\cone'$ and the
multiplication by numbers will also be denoted by $\ps$ and $\pa$
respectively. In the following result, we list the properties that a
coherent map should possess. The mapping that assesses the risk of a
deterministic portfolio in the new cone $\cone'$ will be
$h(f(\cdot))$.  Recall that $\FF$ denotes the family of possible
values of the function $f$.

\begin{proposition}
  \label{risk-transform}
  A map $h:\cone\mapsto\cone'$ is a risk transformation if it is
  \begin{itemize}
  \item[(i)] non-decreasing, i.e. $h(a)\preccurlyeq h(b)$ if
    $a\preceq b$;
  \item[(ii)] linear on $\FF$, i.e. $h(a\ps b)=h(a)\ps h(b)$ for all
    $b\in\cone$ and $a\in\FF$.
  \end{itemize}
  Further, $h$ is a coherent map if $h$ is homogeneous, i.e.  $h(t\pa
  a)=t\pa h(a)$ for all $t>0$ and $a\in\cone$ and also satisfies
  \begin{equation}
    \label{eq:hconvex}
    h(a)\ps h(b)\preccurlyeq h(a\ps b)
  \end{equation}
  for all $a,b\in\cone$.
\end{proposition}
\begin{proof}
  Since $\risk$ satisfies \textbf{R1} and $f(y)\in\FF$, we have for
  all $y\in\R^d$
  \begin{displaymath}
    h\big(\risk(X+y)\big)=h\big(f(y)\ps \risk(X)\big)
    =h\big(g(y)\big)\ps h\big(\risk(X)\big)\,,
  \end{displaymath}
  i.e. \textbf{R1} holds. Property \textbf{R2} holds because $h$ is
  non-decreasing.  The homogeneity of $h(\risk(\cdot))$ is evident if
  $h$ is homogeneous.  If $\risk$ is coherent and (\ref{eq:hconvex})
  holds, then
  \begin{displaymath}
    h\big(\risk(X)\big)\ps h\big(\risk(Y)\big)\preccurlyeq
    h\big(\risk(X)\ps \risk(Y)\big)\preccurlyeq h\big(\risk(X+Y)\big)\,.
  \end{displaymath}\qed
\end{proof}

As an immediate consequence of Proposition~\ref{risk-transform} we
deduce that every linear non-decreasing map is coherent. Such maps
between partially ordered vector spaces are called Riesz
homomorphisms, see \cite[Sec.~18]{lux:zaan71}.

\begin{example}[Vector-valued risk measures from set-valued ones]
  \label{ex:sval-vect}
  A particularly important instance of transformations of risks arises
  if $\coneK$ is a family of convex closed subsets of $\R^d$ with
  inclusion order defined in Example~\ref{ex:jouini} and $\cone'$ is
  $\R^d$ with the reversed $\leqK$-order for a proper Euclidean cone
  $K$. The reversing is needed since $y+K\subseteq z+K$ (i.e.
  $y+K\preceq z+K$) if and only if $z\leqK y$.

  The cone $K$ is said to be a \emph{Riesz cone} if $\R^d$ with
  $\leqK$-order is a Riesz space, i.e.  for every $x,y\in\R^d$ their
  supremum is well defined. It follows
  from~\cite[Th.~26.11]{lux:zaan71} that each Riesz cone can be
  represented as $K=\{u\in\R^d:\; Au\in\R_+^d\}$ for a non-singular
  $d\times d$ matrix $A$ with non-negative entries, i.e.
  $K=A^{-1}\R^d_+$. The matrix $A$ can represent possible transfers
  between the assets so that $Y\leqK X$ if and only if $AY$ is
  coordinatewise smaller than $AX$.

  Assume that $K$ is a Riesz cone. Then it is easy to see that
  $\cone'$ is a complete lattice.  Let $h(F)$ denote the supremum of
  $F\subset\R^d$ in $\cone'$ (i.e. the $\leqK$-infimum of $F$).  If
  $\risk$ is a $\cone$-valued risk measure, then $h(\risk(\cdot))$ is
  a vector-valued risk measure.  Indeed, the map $h$ is monotone and
  homogeneous. Since
  \begin{displaymath}
    h(F-y+K)=h(F-y)=h(F)-y=h(F)+h(-y+K)\,,
  \end{displaymath}
  $h$ is linear on $\FF$. Finally, $h$ satisfies (\ref{eq:hconvex}),
  since $x=h(F_1)$ and $y=h(F_2)$ imply that $F_1+F_2\subseteq
  (x_1+x_2)+K$.

  It is also possible to produce vector-valued risk measures from
  set-valued risk measures in the cone $\coneKR$ from
  Example~\ref{ex:ac-sval} if $h$ is chosen to be the supremum in
  $\cone'$ of the complement of $F\in\coneKR$.
\end{example}

\begin{example}[Linear transformations of vector-valued risk measures]
  \label{ex:lt}
  Let $\risk$ be a risk measure on $\Ldi$ with values in $\R^d$ with
  the reversed $\leqK$-order for a Riesz cone $K$. Note that $K$
  generates both the order on $\Ldi$ and on the space of values for
  $\risk$.  Then
  \begin{displaymath}
    \tilde{\risk}(X)=A^{-1}\risk(AX)
  \end{displaymath}
  is a risk measure with values in $\R^d$ with the reversed
  coordinatewise order.
\end{example}

\begin{example}[Scalar risk measures from vector-valued ones]
  \label{ex:vect-scal}
  Let $K$ be a Riesz cone and $\cone=\R^d$ with the
  reversed $\leqK$-order. Define $\cone'=\R$ with the reversed natural
  order. Finally, let $h(a)=\langle a,u\rangle$, where
  $\langle\cdot,\cdot\rangle$ is the scalar product and $u$ belongs to
  the \emph{positive dual cone} to $K$, i.e.  $\langle u,v\rangle\geq
  0$ for all $v\in K$. Clearly $h$ is a coherent map and we obtain
  univariate risk measures as those of Example~\ref{ex:artzner}, but
  now for multivariate portfolios.
\end{example}

\section{Depth-trimmed regions}
\label{sec:depth}

Depth functions assign to a point its degree of centrality with
respect to the distribution of a random vector, see Zuo and
Serfling~\cite{zuo:ser00a}. The higher the depth of a point is,
the more central this point is with respect to the distribution of
the random vector.  Depth-trimmed (or central) regions are sets of
central points associated with a random vector.  Given a depth
function, depth-trimmed regions can be obtained as its level sets.
With a $d$-dimensional random vector $X$ we associate the family
of \emph{depth-trimmed regions}, i.e. sets $\depth(X)$,
$\alpha\in(0,1]$, such that the following properties hold for all
$\alpha\in[0,1]$ and all $X\in\Ldi$:
\begin{description}
\item[\bf D1] $\depth(X+y)=\depth(X)+y$ for all $y\in\R^d$;
\item[\bf D2] $\depth(tX)=t\depth(X)$ for all $t>0$;
\item[\bf D3] $\depth(X)\subseteq\depth[\beta](X)$ if $\alpha\geq\beta$;
\item[\bf D4] $\depth(X)$ is connected and closed.
\end{description}
Note that the addition of $y$ in \textbf{D1} and the multiplication by
$t$ in \textbf{D2} are the conventional translation and the rescaling
of sets in $\R^d$.

These properties are similar to those discussed by Zuo and
Serfling~\cite[Th.~3.1]{zuo:ser00b}.
Additionally, \cite{zuo:ser00b} requires that the depth-trimmed
regions are invariant with respect to linear transformations, i.e.
$\depth(AX)=A\depth(X)$ for any nonsingular matrix $A$.

We will consider two additional properties of depth-trimmed
regions, that, to our knowledge, have not been studied in the
literature so far:
\begin{description}
\item[\bf D5] if $Y\leqK X$ a.s., then $\depth(X)\subseteq\depth(Y)\pmi
  K$, and $0\in \depth(X)\subseteq K$ if $X=0$ a.s.;
\item[\bf D6] $\depth(X+Y)\subseteq \depth(X)\pmi\depth(Y)$.
\end{description}
Observe that depth-trimmed regions are closed subsets of $\R^d$ and
the addition operation in \textbf{D5} and \textbf{D6} is the closed
Minkowski addition. Later on we will see that \textbf{D6} is closely
related to the coherence property of risk measures.

\begin{example}[Halfspace trimming]
  \label{ex:hs-t}
  The \emph{halfspace trimmed regions} are built as the intersection
  of closed halfspaces whose probability is not smaller than a given
  value:
  \begin{align*}
    \HD^\alpha(X)
    =\bigcap\big\{H:\,H\textrm{ closed halfspace with }\Prob{X\in
    H}\geq1-\alpha\big\}\,.
  \end{align*}
  The above definition of the halfspace trimmed regions is taken
  from Mass\'e and Theodorescu~\cite{mas:theod94}.  Alternatively, the
  strict inequality in the definition of $\HD^\alpha$ is replaced by
  the non-strict one, see Rousseeuw and Ruts~\cite{rous:rut99}.
  However the definition of \cite{mas:theod94} leads to a simpler
  relationship between the value at risk and the univariate halfspace
  trimming, see Section~\ref{sec:accept-sets-gere}.

  It is well known that the halfspace trimmed regions satisfy
  \textbf{D1}--\textbf{D4} and they are compact and convex. The new
  property \textbf{D6} does not hold in general; this can be shown in
  the univariate case using examples for which the value at risk does
  not satisfy \textbf{R4}. The monotonicity property \textbf{D5} does
  neither hold in general.

  However, it is possible to build a variant of the halfspace trimmed
  regions satisfying \textbf{D5}. We define the \emph{monotone
    halfspace trimmed regions} as
  \begin{equation}
    \label{eq:halfspace}
    \HD^\alpha_K(X)=\bigcap_{u\in K^*}
    \big\{H_u:\; \Prob{X\in H_u(t)}\geq 1-\alpha\big\}\,,
  \end{equation}
  where $H_u=\{x\in\R^d:\; \langle x,u\rangle \geq 1\}$ denotes a
  halfspace, and $K^*=\{u:\; \langle u,v\rangle \geq 0,\, v\in K\}$ is
  the positive dual cone to $K$.  The monotone halfspace trimmed
  regions satisfy \textbf{D1}--\textbf{D5} and are nonempty for all
  $\alpha\in(0,1]$.
\end{example}

\begin{example}[Zonoid trimming]
  \label{ex:z-t}
  Koshevoy and Mosler~\cite{kos:mos97z} defined \emph{zonoid trimmed
    regions} for an integrable random vector $X$ in $\R^d$ as
  \begin{equation}\label{eq:zonoid}
    \ZD^\alpha(X)=\big\{\E[X l(X)]:\,l:\R^d\mapsto[0,\alpha^{-1}]
    \textrm{ measurable and }\E l(X)=1\big\}\,,
  \end{equation}
  where $\alpha\in(0,1]$.  Properties \textbf{D1}--\textbf{D4}
  together with convexity and boundedness (and thus compactness) are
  already derived in~\cite{kos:mos97z}. The proofs of \textbf{D5} and
  \textbf{D6} do not involve serious technical difficulties.
\end{example}

\begin{example}[Expected convex hull trimming]
  \label{ex:ch-t}
  \emph{Expected convex hull regions} of a random vector $X$ at level
  $n^{-1}$ for $n\geq1$ are defined by Cascos~\cite{cas05} as the
  selection (Aumann) expectation of the convex hull of $n$
  independent copies $X_1,\dots,X_n$ of $X$, see
  \cite[Sec.~2.1]{mo1} for the definition of expectation for random
  sets.  The expected convex hull region can be given implicitly in
  terms of its support function as
  \begin{displaymath}
    h(\CD^{1/n}(X),u)
    =\E\max\{\langle X_1,u\rangle,\langle X_2,u\rangle,\ldots\langle
    X_n,u\rangle\}\quad\textrm{for all } u\in\R^d\,,
  \end{displaymath}
  where $\langle\cdot,\cdot\rangle$ is the scalar product. Note that
  for any $F\subset\R^d$ its \emph{support function} is given by
  $h(F,u)=\sup\{\langle x,u\rangle:\,x\in F\}$ for $u\in\R^d$.  The
  expected convex hull regions satisfy properties
  \textbf{D1}--\textbf{D6} and are compact and convex.
\end{example}

\begin{example}[Integral trimming]
  \label{ex:int-t}
  Let $\sF$ be a family of measurable functions from $\R^d$ into $\R$.
  Cascos and L\'opez-D\'{\i}az~\cite{cas:lop05} defined the family of
  integral trimmed regions as
  \begin{align*}
    \depth_{\sF}(X)
    &=\bigcup_{Y\in \sP_\alpha(X)}
    \bigg\{x\in\R^d:\, \ftr(x)\leq \E\;\ftr(Y)
    \,\textrm{ for all }\ftr\in\sF\bigg\}\\
    &=\bigcup_{Y\in \sP_\alpha(X)}\;\bigcap_{\ftr\in\sF}
    \ftr^{-1}\big((-\infty,\E\;\ftr(Y)]\big)\,,
  \end{align*}
  where $\sP_\alpha(X)$ is defined in
  Section~\ref{sec:worst-conditioning}.  All families of integral
  trimmed regions satisfy \textbf{D3}.  Other properties of the
  integral trimmed regions heavily depend on their generating family
  of functions. For instance, if for any $\ftr\in\sF$, $t>0$ and
  $z\in\R^d$, the function $\ftr_{t,z}$ defined as
  $\ftr_{t,z}(x)=\ftr(tx+z)$ belongs to $\sF$, then the integral
  trimmed regions generated by $\sF$ satisfy properties \textbf{D1}
  and \textbf{D2}.

  If $\sF=\{\ftr_{t,z}:\,t>0,z\in\R^d\}$ with continuous and
  $\leqK$-decreasing function $\ftr$, then
  \begin{equation}
    \label{eq:itr-da}
    \depth_{\sF}(X)=\bigcup_{Y\in \sP_\alpha(X)}\;
    \bigcap_{t>0\atop z\in\R^d}\;
    \bigg[\frac{1}{t}\Big(\ftr^{-1}
    \big(\E\;\ftr(tY+z)\big)-z\Big)\pmi K\bigg]\,.
  \end{equation}
\end{example}

Hereafter we will assume that all depth-trimmed regions satisfy
\textbf{D1}--\textbf{D5}.

\section{Risk measures generated by depth-trimmed regions}
\label{sec:accept-sets-gere}

As a motivation for the following, note that for an essentially
bounded random variable $X$, $\alpha\in(0,1]$ and $n\geq1$, we have
\begin{align*}
  \VaR(X)&=-\min\HD^\alpha_{[0,\infty)}(X)\,,\\
  \ES(X)&=-\min\ZD^\alpha(X)\,,\\
  \EM(X)&=-\min\CD^{1/n}(X)\,.
\end{align*}
The following example provides another argument showing relationships
between depth-trimmed regions and risk measures.

\begin{example}[Depth-trimmed regions as set-valued risk measures]
  \label{depth-risk}
  Observe that any depth-trimmed region that satisfies
  \textbf{D1}--\textbf{D5} can be transformed into a set-valued risk
  measure from Definition~\ref{def:1}.  Namely,
  $\risk(X)=\depth(X)\pmi K$ is a risk measure in the cone $\coneK$ of
  closed subsets of $\R^d$ with the addition operation being the
  closed Minkowski addition and the \emph{reversed} inclusion order.
  Because of the reversed order, the function $f$ is given by
  $f(x)=x+K$. However, the obtained risk measure is not coherent even
  if \textbf{D6} holds. 

  In order to construct a coherent risk measure from depth-trimmed
  regions, define $\risk(X)$ to be the closure of the complement to
  $\depth(X)\pmi K$. Then $\risk$ becomes a coherent risk measure in
  the cone $\coneKR$ from Example~\ref{ex:ac-sval} if the
  depth-trimmed region satisfies \textbf{D1}--\textbf{D6}.
\end{example}

In general, a random portfolio $X$ will be acceptable or not depending
on the depth-trimmed region of level $\alpha$ associated with $X$.
Since the depth-trimmed regions are subsets of the space $\R^d$ where
$X$ takes its values, we need to map it into the space $\cone$ where
risk measures take their values. This map is provided by the function
$f$ from Definition~\ref{def:1}. Then
\begin{displaymath}
  \drisk(X)=f(\depth(X)\pmi K)
\end{displaymath}
is a subset of $\cone$.  Recall that the acceptance cone $\AA$ is a
subset of $\cone$ that characterises the acceptable values of the risk
measure, see (\ref{f-acceptance}).

\begin{definition}
  \label{accept-set}
  The \emph{acceptance set at level $\alpha$} associated with the
  depth-trimmed region $\depth(\cdot)$ and function $f$ is defined as
  \begin{equation}
    \label{eq:drisk}
    \arisk=\{X\in\Ldi:\,\drisk(X)\subseteq\AA\}.
  \end{equation}
\end{definition}

\begin{theorem}
  \label{th:accept-sets}
  The acceptance sets associated with depth-trimmed regions satisfy
  the following properties:
  \begin{itemize}
  \item[(i)] $0\in\arisk$ for all $\alpha$;
  \item[(ii)] if $\alpha\geq\beta$, then $\arisk[\beta]\subseteq\arisk$;
  \item[(iii)] if $X\in\arisk$, then $tX\in\arisk$ for all $t>0$;
  \item[(iv)] if $X\in\arisk$ and $f(x)\in\AA$, then $x+X\in\arisk$;
  \item[(v)] if $Y\in\arisk$ and $Y\leqK X$ a.s., then $X\in\arisk$;
  \item[(vi)] if $X,Y\in\arisk$ and \textbf{D6} holds, then
    $X+Y\in\arisk$.
  \end{itemize}
\end{theorem}
\begin{proof}
  \noindent \textsl{(i)} By
  \textbf{D5}, $\drisk(0)= f(K)\subseteq\AA$, i.e.
  $0\in\arisk$ for all $\alpha$.

  \noindent \textsl{(ii)} By \textbf{D3},
  $\depth(X)\subseteq\depth[\beta](X)$ whenever $\alpha\geq\beta$.
  Thus $\drisk(X)\subseteq\drisk[\beta](X)$ and
  $\arisk[\beta]\subseteq\arisk$ trivially holds.

  \noindent \textsl{(iii)} By \textbf{D2} and the homogeneity of $f$,
  we have $\drisk(tX)=t\pa \drisk(X)$ for all $t>0$. Since $\AA$ is a
  cone, $\drisk(tX)\subseteq\AA$ if and only if
  $\drisk(X)\subseteq\AA$.

  \noindent \textsl{(iv)} Let $f(x)\in\AA$. By
  \textbf{D1}, we have $\drisk(X+x)=f(\depth(X)+x)$ and by
  (\ref{eq:f}), we have $\drisk(X+x)=\drisk(X)\ps f(x)\subseteq\AA$
  because $\AA$ is a (convex) cone.  By (\ref{eq:drisk}),
  $X+x\in\arisk$.

  \noindent \textsl{(v)} Note that $f(\depth(Y)\ps K)\subseteq\AA$. By \textbf{D5},
  $\depth(X)\ps K\subseteq\depth(Y)\ps K$ and thus $f(\depth(X)\ps K)\subseteq\AA$.

  \noindent \textsl{(vi)} By (\ref{eq:f}) and \textbf{D6},
  \begin{displaymath}
    \drisk(X+Y)\subseteq
    f(\depth(X)\pmi \depth(Y)\ps K)=\drisk(X)\ps
    \drisk(Y)\subseteq\AA\,.
  \end{displaymath}
  Finally, the fact that $\AA$ is a convex cone yields that $X+Y\in\arisk$.
  \qed
\end{proof}

Similarly to the construction used in
Section~\ref{sec:acceptance-sets}, we measure the risk of a portfolio
$X$ in terms of the collection of deterministic portfolios $x$ that
\emph{cancel} the risk induced by $X$ and make $X+x$ acceptable.

\begin{definition}
  \label{def:srisk}
  The risk measure induced by a family of depth-trimmed regions
  $\depth$ at level $\alpha$ is given by
  \begin{equation}
    \label{eq:srisk}
    \srisk(X)=\bigvee \{f(y):\; f(\depth(X-y)\ps K)\subseteq\AA\,, \;
    y\in\R^d\}\,.
  \end{equation}
\end{definition}

By \textbf{D1}, $\srisk(X)$ can be given alternatively in terms of the
acceptance set at level $\alpha$ as
\begin{equation}
  \label{eq:srisk-1}
  \srisk(X)=\bigvee \{f(y):\; X-y\in\arisk\,, \; y\in\R^d\}\,.
\end{equation}

\begin{theorem}
  \label{thr:srisk}
  Assume that $\FF$ is sup-generating. Then the mapping $\srisk(X)$
  satisfies
  \begin{equation}
    \label{eq:sbw}
    \srisk(X)=\bigwedge \drisk(X)\,,
  \end{equation}
  and so becomes a homogeneous risk measure associated with $f$.  If
  the family of depth trimmed regions satisfies \textbf{D6}, then
  $\srisk(X)$ is a coherent risk measure such that
  $\srisk(X)\succeq\srisk[\beta](X)$ for $\alpha\geq\beta$.
\end{theorem}
\begin{proof}
  The linearity of $f$ and (\ref{eq:srisk}) imply that
  \begin{align*}
    \srisk(X)&=\bigvee\{f(y):\; f(\depth(X)\ps K)\subseteq\AA\ps f(y)\}\\
    &=\bigvee\{f(y):\;c\in\AA\ps f(y)\;
    \textrm{ for all }\;c\in\drisk(X)\}\\
    &=\bigvee\{f(y):\; f(y)\preceq c\;
    \text{ for all }\; c\in\drisk(X)\}\\
    &=\bigvee\{a\in\FF:\; a\preceq \bigwedge\drisk(X)\}\,,
  \end{align*}
  so that (\ref{eq:sbw}) follows from the sup-generating property
  (\ref{eq:sg}).

  If $X=x$ a.s., then $\srisk(X)=\bigwedge f(x+\depth(0)\ps K)=f(x)$,
  since $\bigwedge f(\depth(0)\ps K)=\neutral$ by \textbf{D5} and $f$
  is non-decreasing.  By (\ref{eq:f}) and \textbf{D1} we deduce that
  $\srisk(X)\ps f(y)=\srisk(X+y)$, so \textbf{R1} holds.

  If $Y\leqK X$ a.s., then $\depth(X)\subseteq\depth(Y)\ps K$ by
  \textbf{D5}. Thus \textbf{R2} holds, since
  \begin{displaymath}
    \srisk(Y)=\bigwedge f(\depth(Y)\ps K)\preceq\bigwedge f(\depth(X)\ps K)=\srisk(X)\,.
  \end{displaymath}
  Property \textbf{R3} follows directly from \textbf{D2}, the fact
  that $K$ is a cone and the homogeneity of $f$.  If \textbf{D6}
  holds, then $\depth(X+Y)\ps K\subseteq\left(\depth(X)\ps
    K\right)\ps\left(\depth(Y)\ps K\right)$
  \begin{displaymath}
    \srisk(X+Y)=\bigwedge\drisk(X+Y)\succeq\bigwedge\drisk(X)\ps\bigwedge\drisk(Y)
    =\srisk(X)\ps\srisk(X)\,,
  \end{displaymath}
  i.e. \textbf{R4} holds. Finally, the ordering of the risks with
  respect to $\alpha$ follows from \textbf{D3}.\qed
\end{proof}


Now we describe a dual construction, based on rejection sets, of
set-valued risk measures associated with depth-trimmed regions.  The
rejection set at level $\alpha$ associated with $\depth(\cdot)$ is
given by
\begin{displaymath}
  \ariskR=\{X\in\Ldi:\; \drisk(X)\cap\AAR\neq\emptyset\}
  =\{X\in\Ldi:\; \depth(X)\cap(-K)\neq\emptyset\}\,.
\end{displaymath}
Assuming that $\FF$ is inf-generating, by (\ref{eq:r-ar}) we have
\begin{displaymath}
  \sriskR(X)=\bigwedge\risk_{\ariskR}(X)\,,
\end{displaymath}
where $\sriskR(X)$ is also given by (\ref{eq:srisk-1}) with $\arisk$
replaced by $\ariskR$.  It is possible to reproduce
Theorem~\ref{thr:srisk} in this dual framework and obtain that
\begin{displaymath}
  \sriskR(X)=\bigwedge\drisk(X)\,.
\end{displaymath}
Further, $\sriskR$ is a homogeneous risk measure which is also
coherent if \textbf{D6} holds.

\begin{example}[Set-valued risk measures from depth-trimmed regions]
  \label{ex:sva-d}
  In the setting of Example~\ref{ex:mvar-cones} $f(x)=-x+K$, so that
  Theorem~\ref{thr:srisk} implies that
  \begin{align*}
    \srisk(X)&=\bigcap_{x\in\depth(X)} (-x+K)
    =\bigcap_{x\in\depth(X)} \{z\in\R^d:\; -x\leqK z\}\\
    &=\{z\in\R^d: z+\depth(X)\subseteq K\}
    =\{z\in\R^d:\; \depth(X)\subseteq (-z+K)\}\,.
  \end{align*}
  If $K$ is a Riesz cone, then there exists the infimum of $\depth(X)$
  with respect to the $\leqK$-order (denoted as $\wedge_K\depth(X)$),
  so that
  \begin{equation}
    \label{eq:new-sva}
    \srisk(X)=\bigcap_{x\in\depth(X)} \{z\in\R^d:\; -z\leqK x\}
    =-\wedge_K\depth(X)+K\,.
  \end{equation}
  Therefore, risk measures generated by depth-trimmed regions using
  the acceptance cone are not particularly interesting, since they are
  essentially vector-valued. In Example~(\ref{eq:mrep}) it will be
  shown that vector-valued risk measures are necessarily marginalised,
  i.e. they appear from the scheme of Example~\ref{ex:vector}.

  However, the rejection construction produces more interesting
  set-valued risk measures. Namely, in the setting of
  Example~\ref{ex:mvar-cones} with $f(x)=x+\KR$, the corresponding
  risk measure $\sriskR(X)$ is the closure to the complement of
  $\depth(X)\pmi K$. The obtained risk measure takes values in the
  cone $\coneKR$.
\end{example}

\section{Basic risk measures associated with depth-trimmed regions}
\label{sec:basic-risk-measures}

Let us now specialise the constructions from
Section~\ref{sec:accept-sets-gere} for $X=(X_1,X_2,\dots,X_d)\in\Ldi$
and several basic definitions of depth-trimmed regions and set-valued
risk measures with values either in $\coneK$ with a Riesz cone $K$.
Recall that set-valued risk measures with values in $\coneK$ can be
represented as $x+K$ for some $x\in\R^d$, i.e. are effectively
vector-valued.  Similar constructions are possible for
$\coneKR$-valued risk measures.  In this case the $\leqK$-infimum of
the complement to $\coneKR$-valued risk measures also yields a
vector-valued risk measure, see Example~\ref{ex:sval-vect}.



\paragraph{Risk measures generated by monotone halfspace trimming.}
The monotone halfspace trimming induces a homogeneous risk
measure, i.e. \textbf{R3} holds. This set-valued risk measure is
given by $\srisk(X_1)=\big[\VaR(X_1),+\infty\big)$ in the
univariate case. In general,
\begin{align*}
  \srisk(X_1,X_1,\dots,X_1)
  =&\big(\VaR(X_1),\VaR(X_1),\dots,\VaR(X_1)\big)+K\,,\\
  \srisk(X)
  \supseteq&\big(\VaR(X_1),\VaR(X_2),\dots,\VaR(X_d)\big)+K\,.
\end{align*}

\paragraph{Risk measures generated by zonoid trimming.}
The zonoid trimming induces coherent risk measures.  Then
$\srisk(X_1)=\big[\ES(X_1),+\infty\big)$ and in the multivariate
setting
\begin{align*}
  \srisk(X_1,X_1,\dots,X_1)
  &=\big(\ES(X_1),\ES(X_1),\dots,\ES(X_1)\big)+K\,,\\
  \srisk(X)
  &\supseteq\big(\ES(X_1),\ES(X_2),\dots,\ES(X_d)\big)+K\,,
\end{align*}
where the latter inclusion turns into the equality if
$K=\R^d_+\,$. If $K=A^{-1}\R_+^d$ for a matrix $A$, then
\begin{equation}
  \label{eq:es-mvar}
  \srisk(X)=A^{-1}\big(\ES((AX)_1),\ES((AX)_2),\dots,\ES((AX)_d)\big)+K\,,
\end{equation}
where $(AX)_i$ stands for the $i$-th coordinate of $AX$. In
particular, (\ref{eq:es-mvar}) implies that the marginalised expected
shortfall (as in Example~\ref{ex:vector}) of $AX$ is coordinatewise
smaller than $A$ applied to the marginalised expected shortfall of
$X$.

\paragraph{Risk measures generated by expected convex hull trimming.}
The expected convex hull trimming induces coherent risk measures.  Then
$\sriskn(X_1)=\big[\EM(X_1),+\infty\big)$ and
\begin{align*}
  \sriskn(X_1,X_1,\dots,X_1)
  &=\big(\EM(X_1),\EM(X_1),\dots,\EM(X_1)\big)+K\,,\\
  \sriskn(X)
  &\supseteq\big(\EM(X_1),\EM(X_2),\dots,\EM(X_d)\big)+K
\end{align*}
with the equality if $K=\R^d_+\,$. If $K=A^{-1}\R_+^d$ for a
matrix $A$, then (\ref{eq:es-mvar}) also holds for the expected
minimum instead of the expected shortfall.

Note that in all three examples described above we have 
\begin{displaymath}
  \sriskR(X)\supseteq (\rho(X_1),\dots,\rho(X_d))+\KR\,,
\end{displaymath}
where $\rho$ stands for $\VaR$, $\ES$ or $\EM$. 


\paragraph{Integral trimmed risk measures.}
The integral trimmed regions generate new multivariate risk measures.
Consider the cone $\coneK$ from Example~\ref{ex:jouini} and
$f(x)=-x+K$.  Let $\sF=\{\ftr(tx+z):\,t>0,\,z\in\R^d\}$, where $\ftr$
is continuous and $\leqK$-decreasing for a proper Riesz cone $K$.
Since $\depth(X)$ is the union of
\begin{displaymath}
  \depth[1](Y)=\bigcap_{t>0\atop z\in\R^d}\left[
  \frac{1}{t}\Big(\ftr^{-1}\big(\E\;\ftr(tY+z)\big)-z\Big)\pmi
      K\right]
\end{displaymath}
for $Y\in\sP_\alpha(X)$, we obtain
\begin{displaymath}
  \srisk(X)=\bigcap_{Y\in\sP_\alpha(X)} \srisk[1](Y)\,,
\end{displaymath}
so that $\srisk(X)$ appears from the worst conditioning
construction applied to the risk measure $\srisk[1]$. Furthermore,
(\ref{eq:new-sva}) yields that $\srisk[1](X)=x+K$, where
\begin{align}
  x=-\wedge_K\depth[1](X)
  &=-\wedge_K \left(\bigcap_{t>0\atop z\in\R^d}\; \left[
      \frac{1}{t}\Big(\ftr^{-1}\big(\E\;\ftr(tY+z)\big)-z\Big)\pmi
      K\right]\right)\nonumber\\
  &=-\bigvee_{t>0\atop z\in\R^d}\wedge_K\left[
    \frac{1}{t}\Big(\ftr^{-1}\big(\E\;\ftr(tY+z)\big)-z\Big)\right]
  \nonumber\\
  &=\bigwedge_{t>0\atop z\in\R^d}\left(-\wedge_K\left[
      \frac{1}{t}\Big(\ftr^{-1}\big(\E\;\ftr(tY+z)\big)
      -z\Big)\right]\right)\,\label{eq:srk1}.
\end{align}
This risk measure satisfies \textbf{R1}--\textbf{R3} and results from
the homogenisation construction (\ref{eq:hc}) and (\ref{eq:htc})
applied to the set-valued risk measure generated by the integral
trimmed regions whose generating family is $\sF=\{\ftr\}$,
\begin{equation}
  \label{eq:srk2}
  \risk(X)=-\wedge_K \ftr^{-1}\big(\E\;\ftr(X)\big)+K\,.
\end{equation}
Notice that this homogenisation preserves \textbf{R2}, but not
necessarily \textbf{R4}. The idea of constructing scalar risk measures
using real-valued functions of vector portfolios appears also in
\cite{bur:rues06}.  Alternatively, it is possible to take infimum in
over $t>0$ or over $z\in\R^d$ only, which results in a risk measure
that satisfies \textbf{R3} or \textbf{R1} respectively.

\begin{example}
  The function $\ftr(t)=e^{-t/\gamma}$ yields the risk measure
  $\risk(X)=\gamma\log(\E e^{-X/\gamma})$ by (\ref{eq:srk2}) in
  $\cone=\R$ with the reversed order and $f(x)=-x$. The properties
  \textbf{R1} and \textbf{R2} evidently holds, while (\ref{eq:crm})
  follows from the H\"older inequality, i.e. $\risk$ is a
  \emph{convex} risk measure, which does not satisfy \textbf{R3}.
  Since \textbf{R1} already holds, there is no need to take infimum
  over $z\in\R^d$ in (\ref{eq:srk1}). The corresponding convex risk
  measure is called the \emph{entropic risk measure} with $\gamma$
  being the risk tolerance coefficient.

  If we attempt to produce a homogeneous (and thereupon coherent) risk
  measure from $\risk$, we need to apply (\ref{eq:hc}), which in view
  of the reversed order on the real line turns into
  \begin{displaymath}
    \risk_{\mathrm{h}}(X)=\sup_{t>0}\; t^{-1}\, \risk(tX)
    =\sup_{t>0}\; t^{-1}\,\log (\E e^{-tX})
    =\sup_{t>0}\; \log ((\E Y^t)^{1/t})
  \end{displaymath}
  for $Y=e^{-X}$. Since $(\E Y^t)^{1/t}$ is an increasing function of
  $t>0$, we have
  \begin{displaymath}
    \risk_{\mathrm{h}}(X)=\lim_{t\to\infty}
    t^{-1}\log (\E e^{-tX})\,.
  \end{displaymath}
  It is easy to see that the limit equals $(-\essinf X)$, so a
  coherent variant of $\risk$ is not particularly interesting.

\end{example}

\section{Depth-trimmed regions generated by risk measures}
\label{sec:depth-trimm-regi}

Consider a family of homogeneous risk measures $\risk_\alpha$ for
$\alpha\in(0,1]$ such that
\begin{equation}
  \label{eq:rmon}
  \risk_\alpha\succeq\risk_\beta\,,\quad \alpha\geq\beta\,,
\end{equation}
which are associated with function $f$ according to \textbf{R1}. For
instance, such family of risk measures can be produced using the worst
conditioning construction from Section~\ref{sec:worst-conditioning}.

\begin{definition}
  \label{def:back-depth}
  The \emph{depth-trimmed regions generated by the family of risk
    measures} are defined as
  \begin{displaymath}
    \depth(X)=\{y\in\R^d:\; \risk_\alpha(X-y)\preceq\neutral\}\,.
  \end{displaymath}
\end{definition}

By \textbf{R1}, the depth-trimmed regions generated by a family of
risk measures are alternatively given by
\begin{equation}
  \label{eq:depth:risk}
  \depth(X)=\{y\in\R^d:\;\risk_\alpha(X)\preceq f(y)\}\,.
\end{equation}

\begin{theorem}
  \label{thr:depth-risk}
  The depth-trimmed regions generated by a family of risk measures satisfy
  \begin{itemize}
  \item[(i)] properties \textbf{D1}, \textbf{D2}, \textbf{D3} and
    \textbf{D5};
  \item[(ii)] are convex if the risk measure is convex;
  \item[(iii)] are closed if $f$ is upper semicontinuous, i.e.
    $\{x\in\R^d:\; a\preceq f(x)\}$ is closed in $\R^d$ for every
    $a\in\cone$.
  \end{itemize}
\end{theorem}
\begin{proof}
  \noindent\textsl{(i)} Properties \textbf{D1} and \textbf{D2}
  trivially hold by \textbf{R1} and \textbf{R3} respectively.  The
  nesting property \textbf{D3} of depth-trimmed regions is a
  consequence of (\ref{eq:rmon}). We will show that \textbf{D5}
  follows from \textbf{R2}.  If $Y\leqK X$ a.s., then \textbf{R2}
  yields that $\risk_\alpha(Y)\preceq\risk_\alpha(X)$.  Then
  \begin{displaymath}
    \{y\in\R^d:\,\risk_\alpha(Y)\preceq f(y)\}
    \supseteq\{y\in\R^d:\,\risk_\alpha(X)\preceq f(y)\}
  \end{displaymath}
  and by~(\ref{eq:depth:risk}) we have $\depth(Y)\supseteq\depth(X)$ and
  finally $\depth(X)\subseteq\depth(Y)\pmi K$, since $0\in K$.

  \noindent\textsl{(ii)} Given $y,z\in\depth(X)$ and $t\in[0,1]$,
  \begin{displaymath}
    \risk_\alpha(X)\preceq t\pa f(y)\ps (1-t)\pa f(z)=f(ty+(1-t)z)
  \end{displaymath}
  and finally $ty+(1-t)z\in\depth(X)$.

  \noindent\textsl{(iii)} If $f$ is upper semicontinuous, the set
  $\{y\in\R^d:\;\risk_\alpha(X)\preceq f(y)\}$ is closed. \qed
\end{proof}

Under mild conditions, it is possible to recover a risk measure
from the depth-trimmed regions generated by it. If $\FF$ is
sup-generating and inf-generating, the original risk measure is
the risk measure induced by the family of depth-trimmed regions
that it generated. Theorem~\ref{thr:srisk} and
equation~(\ref{eq:depth:risk}) yield that
\begin{displaymath}
  \srisk(X)=\bigwedge \drisk(X)
  =\bigwedge\{f(y):\;\risk_\alpha(X)\preceq f(y)\}
  =\risk_\alpha(X)\,.
\end{displaymath}
Notice that if the construction based on rejection sets is used,
see (\ref{eq:r-ar}), the first equality in the left holds when
$\FF$ is inf-generating, so the sup-generating assumption on $\FF$
can be dropped and we still have
\begin{displaymath}
  \sriskR(X)=\risk_\alpha(X)\,.
\end{displaymath}

\begin{example}[Expected convex hull trimming revisited]
  \label{ex:ectr-rev}
  The expected minimum can be formulated as a \emph{spectral risk
    measure}, see~\cite{acer02}, as
  \begin{equation}
    \label{eq:minimum}
    \EM(X)=-\int_0^1 n(1-t)^{n-1}F_X^{-1}(t)dt\,,\quad n\geq1\,.
  \end{equation}
  For any $\alpha\in(0,1]$, define $\EM[\alpha](X)$ substituting $n$
  by $\alpha^{-1}$ in~(\ref{eq:minimum}). The risk measures
  $\EM[\alpha]$ generates a family of depth-trimmed regions for
  $X\in\Ldi[1]$ with a continuous parameter $\alpha\in(0,1]$.
  Applying Definition~\ref{def:back-depth}, we obtain
  $\depth(X)=[-\EM[\alpha](X),+\infty)$.

  In contrast to risk measures, depth-trimmed regions treat all
  directions in the same way, so that the regions $\depth$ must be
  slightly modified so that they yield the expected convex hull
  trimmed regions for $\alpha=1/n$.  Define
  \begin{align*}
    \CD^\alpha(X)&=\depth(X)\;\cap\;\big(-\depth(-X)\big)\\
    &=[-\EM[\alpha](X)\,,\,\EM[\alpha](-X)]\\
    &=\bigg[\alpha^{-1}\int_0^1(1-t)^{\alpha^{-1}-1}F_X^{-1}(t)dt\,,
    \,\alpha^{-1}\int_0^1 t^{\alpha^{-1}-1}F_X^{-1}(t)dt\bigg]\,.
  \end{align*}
  In this formulation, we can assume that the parameter $\alpha$ takes
  any value in $(0,1]$ and thus, we obtain an extension of the
  univariate expected convex hull trimmed regions.
\end{example}




\section{Duality results}
\label{sec:duality-results}

The dual space to $\Ldi$ is the family of finitely additive bounded
\emph{vector} measures $\mu=(\mu_1,\dots,\mu_d)$ on the underlying
probability space $(\Omega,\salg)$, which act on $X\in\Ldi$ as
$\E_\mu(X)=\sum_{i=1}^d \int X_id\mu_i$.  The polar set to the cone of
acceptable portfolios to a coherent risk measure $\risk$ can be
written as
\begin{displaymath}
  \sA^*=\bigcap_{X\in\sA}\{\mu:\,\E_\mu(X)\geq 0\}\,.
\end{displaymath}
As in \cite{jouin:med:touz04}, we can apply the bipolar theorem to
show that
\begin{equation}
  \label{eq:sa-rep}
  \sA=\bigcap_{\mu\in\sA^*}\{X:\,\E_\mu(X)\geq 0\}\,.
\end{equation}
For each $u\in K$ and measurable $\Omega'\subset\Omega$, the random
element $u\Ind_{\Omega'}$ belongs to $\sA$. Therefore, every
$\mu\in\sA^*$ satisfies $\sum \mu_i(\Omega')u_i\geq0$ for every $u\in
K$. Thus, the values of any $\mu\in\sA^*$ belong to the positive dual
cone to $K$.

Assume that $\FF$ is sup-generating. Proposition~\ref{prop:r-a}
implies that
\begin{displaymath}
  \risk(X)=\bigvee\{f(y):\; X-y\in \sA\}\,.
\end{displaymath}
It follows from (\ref{eq:sa-rep}) that
\begin{align*}
  \risk(X)&=\bigvee
  \{f(y):\; y\in\R^d,\;
  \E_\mu(X)\geq \langle \mu,y\rangle\; \text{for all}\;
  \mu\in\sA^*\}\\
  &=\bigvee f\Big(\bigcap_{\mu\in\sA^*}
  \{y\in\R^d:\; \langle \mu,y\rangle\leq \E_\mu(X)\}\Big)\,,
\end{align*}
where $\langle \mu,x\rangle=\sum_{i=1}^d \mu_i(\Omega)x_i$.

For instance, a set-valued risk measure with values in the cone of
convex closed sets in $\R^d$ with the inclusion order (so that
$f(-x)=x+K$) can be represented as
\begin{displaymath}
  \risk(X)=\bigcap_{\mu\in\sA^*} \{x\in\R^d:\;
  \langle \mu,x\rangle \geq \E_\mu(-X)\}\,,
\end{displaymath}
where $\sA^*$ is a set of finitely additive bounded vector measures
with values in $K$. Note that there is no need to add $K$ to the
right-hand side, since $\langle \mu,x+z\rangle\geq \langle
\mu,x\rangle$ in case $\mu$ takes values from the positive dual to
$K$. By applying to this set-valued risk measure $\risk$ the worst
conditioning construction, we obtain
\begin{align*}
  \risk_\alpha(X)&=\bigwedge_{Y\in\sP_\alpha(X)}\;
  \bigcap_{\mu\in\sA^*} \{x\in\R^d:\; \langle \mu,x\rangle
  \geq \E_\mu(-Y)\}\\
  &=\bigcap_{\mu\in\sA^*}
  \left\{x\in\R^d:\; \langle \mu,x\rangle\geq
    (-\E_\mu)_\alpha(X)\right\}\,,
\end{align*}
where
\begin{displaymath}
  (-\E_\mu)_\alpha(X)=(-\E_{\mu_1})_\alpha X_1+\cdots+
  (-\E_{\mu_d})_\alpha X_d\,.
\end{displaymath}
Thus $\risk_\alpha$ also admits the dual representation, where instead
of the expectation $\E_\mu(-X)$ we take the expected shortfall of $X$
with respect to the measure $\mu$.  Definition~\ref{def:back-depth}
yields then a dual representation for the family of depth-trimmed
regions.

If the risk measure satisfies the Fatou property, then all measures
from $\sA^*$ can be chosen to be $\sigma$-additive. Recall that the
Fatou property means that the risk measure is lower semicontinuous in
probability, i.e. the lower limit (which for set-valued risk measures
is understood in the Painlev\'e-Kuratowski sense \cite[Def.~B.5]{mo1})
of $\risk(X_k)$ is not greater than $\risk(X)$ if $X_k$ converges in
probability to $X$.

\begin{example}[Vector-valued coherent risk measures]
  \label{ex:degen}
  Let $\risk$ be a coherent risk measure with values in $\cone=\R^d$
  with the reversed $\leqK$-order with $K=\R^d_+$ and $f(x)=-x$.  Then
  \begin{displaymath}
    \risk(X)=\bigvee \bigcap_{\mu\in\sA^*}
    Y(\mu)\,,
  \end{displaymath}
  where all $\mu\in\sA^*$ have non-negative coordinates and 
  \begin{displaymath}
    Y(\mu)=\{-y:\; \langle \mu,y\rangle\leq \E_\mu(X),\; y\in\R^d\}
    =\{y:\; \langle \mu,y\rangle\geq -\E_\mu(X),\; y\in\R^d\}\,.
  \end{displaymath}
  For instance, the first coordinate of $\risk(X)$ is obtained as the
  infimum of the projection of $\cap_{\mu\in\sA^*} Y(\mu)$ onto the
  first coordinate. If $\mu_1(\Omega)>0$ and
  $\mu_2(\Omega)=\cdots=\mu_d(\Omega)=0$, then
  $Y(\mu)=[y_1,\infty)\times\R\cdots\times\R$ for some $y_1$. If $x^*$
  is the essential lower bound of $X$ with respect to $\leqK$-order,
  then $(y_1,-x_2^*,\dots,-x_d^*)$ belongs to $Y(\mu)$ for all
  $\mu\in\sA^*$. Thus, the first coordinate of $\risk(X)$ is given by
  the infimum $-\mu_1(\Omega)^{-1}\E_{\mu_1} X_1$ over all measures
  $(\mu_1,0,\dots,0)$ that belong to $\sA^*$.
  Thus,
  \begin{equation}
    \label{eq:mrep}
    \risk(X)=
    \Big(-\inf_{\mu_1\in\sA^*_1} \E_{\mu_1}(X_1),\dots,-\inf_{\mu_d\in\sA^*_d}
    \E_{\mu_d}(X_d)\Big)\,,
  \end{equation}
  where $\sA^*_i$ is the family of normalised measures $\mu_i$ such
  that $(0,\dots,0,\mu_i,0,\dots,0)\in\sA^*$, $i=1,\dots,d$. The
  individual infima in (\ref{eq:mrep}) are risk measures themselves.
  Therefore, $\risk(X)$ can be represented as the vector composed of
  coherent risk measures of the marginals of $X$.
  
  A similar argument is applicable for the risk measure
  $A^{-1}\risk(AX)$ if $K$ is a general Riesz cone with
  $K=A^{-1}\R_+^d$. Then $A\risk(X)$ can be represented as the vector
  composed of risk measures calculated for the coordinates of $AX$,
  cf. (\ref{eq:es-mvar}).
\end{example}


\section{Conclusions}
\label{sec:conclusions}

It is likely that other results from the morphological theory of
lattices \cite{hei94} have applications in the framework of risk
measures. In particular, it would be interesting to find a financial
interpretation for dilation mappings that commute with supremum,
erosions that commute with infimum, and pairs of erosions and
dilations that are called adjunctions.

It is possible to consider a variant of \textbf{R2} where $Y\leq X$ is
understood with respect to any other chosen order on $\Ldi$. The
consistency issues for scalar risk measures for vector portfolios are
investigated in \cite{bur:rues06} and in  \cite{baeuer:muel06} for the
one-dimensional case.




\section*{Acknowledgements}
\label{sec:acknowledgements}

Ignacio Cascos acknowledges the hospitality of the Department of
Mathematical Statistics and Actuarial Science of the University of
Berne. The encouragement of the editors has led to a substantial
improvement of the manuscript.

\bibliographystyle{spmpsci}      


\newcommand{\noopsort}[1]{} \newcommand{\printfirst}[2]{#1}
  \newcommand{\singleletter}[1]{#1} \newcommand{\switchargs}[2]{#2#1}

\end{document}